# POSITIVITY FOR CERTAIN WEYL GROUP REPRESENTATIONS

G. Lusztig

Abstract. Let W be a Weyl group. We define the notion of positivity of a W-module in terms of the corresponding module over the asymptotic Iwahori-Hecke algebra. Assuming that $W$ is of classical type, we describe a large collection of $W$-modules which are positive. This includes the special representations of $W$, the $W$-modules carried by the left cells of $W$ and also some $W$-modules which interpolate between these two kinds of modules.

## Introduction

**0.1.** Let $\mathcal{G}$ be a finite group acting on a finite set $X$. Then $\mathcal{G}$ acts diagonally on $X^2 := X \times X$ and we can consider the equivariant $K$-theory group $K_{\mathcal{G}}(X^2)$. In [L87a] a convolution operation on $K_{\mathcal{G}}(X^2)$ is defined; it makes $K_{\mathcal{G}}(X^2)$ an (associative) ring with 1. Let $\underline{K}_{\mathcal{G}}(X^2) = \mathbf{C} \otimes K_{\mathcal{G}}(X^2)$; this is an associative algebra over $\mathbf{C}$. Now $\underline{K}_{\mathcal{G}}(X^2)$ has as $\mathbf{C}$-basis

(a) the set of irreducible $\mathcal{G}$-vector bundles on $X^2$ (up to isomorphism).

If $M$ is a $\underline{K}_{\mathcal{G}}(X^2)$-module of finite dimension over $\mathbf{C}$, we say that $M$ is *positive* if there exists a $\mathbf{C}$-basis of $M$ such that the action of any of the elements in (a) in this basis of $M$ is given by a matrix with all entries in $\mathbf{R}_{\geq 0}$.

**0.2.** The main result of this paper is a definition in the case where $\mathcal{G} = E$ (a finite dimensional $\mathbf{Z}/2$-vector space) of a family of $\underline{K}_{\mathcal{G}}(X^2)$-modules that are positive (see 4.7(a)).

**0.3.** We now assume that $W$ is a Weyl group. In §7 we define the notion of positive representation of $W$ (by reduction to the analogous notion for the modules over an algebra as in 0.1). In the case where $W$ is of classical type we define a large collection of $W$-modules which are positive (using a reduction to 4.7(a)). This family includes the special representations of $W$, the $W$-modules carried by the left cells of $W$ and also some $W$-modules which interpolate between these two kinds of modules.

**0.4.** If $\Pi$ is a property (such as an equality or an inclusion) we write $\delta_\Pi = 1$ if $\Pi$ is true and $\delta_\Pi = 0$ if $\Pi$ is not true.

## 1. The algebra $\mathcal{F}$

**1.1.** Let $E$ be a finite dimensional $\mathbf{Z}/2$-vector space. Let $\mathcal{L}$ be the set of vector subspaces of $E$.

Let $X$ be a finite set with a given $E$-action $e : x \mapsto e+x$. We have $X = \sqcup_{A\in\mathcal{L}} X_A$ where $X_A = \{x \in X; \{e \in E; e + x = x\} = A\}$. Let $\bar{X}$ be the set of $E$-orbits in $X$. We have $\bar{X} = \sqcup_{A\in\mathcal{L}} \bar{X}_A$ where $\bar{X}_A$ is the set of $E$-orbits contained in $X_A$. Let $\mathcal{L}_X = \{A \in \mathcal{L}; X_A \neq \emptyset\}$.

Now $E$ acts on $X^2 = X \times X$ by $e : (x, x') \mapsto (e + x, e + x')$. We have $X^2 = \sqcup_{(A,B)\in\mathcal{L}\times\mathcal{L}} X^2_{AB}$ where $X^2_{AB} = X_A \times X_B$. Let $\overline{X^2}$ be the set of $E$-orbits in $X^2$. We have $\overline{X^2} = \sqcup_{(A,B)\in\mathcal{L}\times\mathcal{L}} \overline{X^2}_{AB}$ where $\overline{X^2}_{AB}$ is the set of $E$-orbits contained in $X^2_{AB}$. For any $\Sigma \subset X^2$ we denote by $p_1\Sigma$ (resp. $p_2\Sigma$) the image of $\Sigma$ under the first (resp. second) projection $X^2 \to X$.

**1.2.** For any $\mathbf{Z}/2$-vector space $V$ let $V^*$ be the vector space of linear forms $V \to \mathbf{Z}/2$.

For $A, B$ in $\mathcal{L}$ we consider the $E$-action on $X^2_{AB} \times (A \cap B)^*$ given by

$$e : ((x, y), \epsilon) \mapsto ((e + x, e + y), \epsilon);$$

let $\mathcal{F}_{AB}$ be the $\mathbf{C}$-vector space consisting of all functions $f : X^2_{AB} \times (A\cap B)^* \to \mathbf{C}$ that are constant on the orbits of this $E$-action. Let

$$\mathcal{F} = \oplus_{A,B \text{ in } \mathcal{L}} \mathcal{F}_{AB},$$

a $\mathbf{C}$-vector space.

For $A, B, C$ in $\mathcal{L}$ and $\epsilon_1 \in (A \cap B)^*, \epsilon_2 \in (B \cap C)^*, \epsilon_3 \in (A \cap C)^*$, we write $\bar{\epsilon}_1 + \bar{\epsilon}_2 + \bar{\epsilon}_3 = 0$ if the sum of the restrictions of $\epsilon_1, \epsilon_2, \epsilon_3$ to $A \cap B \cap C$ is zero.

We define a $\mathbf{C}$-bilinear map $\mathcal{F}_{AB} \times \mathcal{F}_{BC} \to \mathcal{F}_{AC}$, $(f_1, f_2) \mapsto f_1 \star f_2$ by

$$(f_1 \star f_2)((x, z), \epsilon) = \sum_{(y,\epsilon_1,\epsilon_2)\in X_B\times(A\cap B)^*\times(B\cap C)^*;\bar{\epsilon}_1+\bar{\epsilon}_2+\bar{\epsilon}=0} f_1((x, y), \epsilon_1)f_2((y, z), \epsilon_2)\frac{|A \cap B \cap C|}{|A \cap C|}$$

for $(x, z) \in X^2_{AC}, \epsilon \in (A \cap C)^*$.

We extends this to a $\mathbf{C}$-bilinear map (product) $\mathcal{F}\times\mathcal{F}\to\mathcal{F}$ by requiring that the product of $\mathcal{F}_{AB}$ and $\mathcal{F}_{B'C}$ is zero if $B\ne B'$. We show:

(a) *For $A,B,C,D$ in $\mathcal{L}$ and for $f_1\in\mathcal{F}_{AB}$, $f_2\in\mathcal{F}_{BC}$, $f_3\in\mathcal{F}_{CD}$ we have $(f_1\star f_2)\star f_3=f_1\star(f_2\star f_3)$.*

Let $(x,u)\in X^2_{AD}$, $\epsilon\in(A\cap D)^*$. We have

$$((f_1\star f_2)\star f_3)((x,u),\epsilon)$$

$$=\sum_{z\in X_C;\epsilon_{12},\epsilon_3}(f_1\star f_2)((x,z),\epsilon_{12})f_3((z,u),\epsilon_3)|A\cap C\cap D|/|A\cap D|$$

$$=\sum_{y\in X_B,z\in X_C;\epsilon_{12},\epsilon_1,\epsilon_2,\epsilon_3}f_1((x,y),\epsilon_1)f_2((y,z)\epsilon_2)f_3((z,u),\epsilon_3)\frac{|A\cap C\cap D}{|A\cap D|}\frac{|A\cap B\cap C}{|A\cap C|}$$

where $(\epsilon_{12},\epsilon_1,\epsilon_2,\epsilon_3)$ runs through the subset

$$\mathfrak{T}\subset(A\cap C)^*\times(A\cap B)^*\times(B\cap C)^*\times(C\cap D)^*$$

defined by the conditions that $\bar\epsilon+\bar\epsilon_{12}+\bar\epsilon_3=0$ and $\bar\epsilon_{12}+\bar\epsilon_1+\bar\epsilon_2=0$.

Let $\mathfrak{T}'$ be the set of all

$$(\epsilon_1,\epsilon_2,\epsilon_3)\in(A\cap B)^*\times(B\cap C)^*\times(C\cap D)^*$$

such that the sum of the restrictions of $\epsilon,\epsilon_1,\epsilon_2,\epsilon_3$ to $A\cap B\cap C\cap D$ is zero. We have a well defined map $\mathfrak{T}\to\mathfrak{T}'$, $(\epsilon_{12},\epsilon_1,\epsilon_2,\epsilon_3)\mapsto(\epsilon_1,\epsilon_2,\epsilon_3)$. The fibre of this map at $(\epsilon_1,\epsilon_2,\epsilon_3)$ can be identified with the set of all $\epsilon_{12}\in(A\cap C)^*$ whose restriction to $A\cap C\cap D$ is the sum $\epsilon'_{12}$ of the restrictions of $\epsilon$ and $\epsilon_3$ and whose restriction to $A\cap B\cap C$ is the sum $\epsilon''_{12}$ of the restrictions of $\epsilon_1$ and $\epsilon_2$. From the definition we see that $\epsilon'_{12}=\epsilon''_{12}$ on $(A\cap C\cap D)\cap(A\cap B\cap C)$. Hence our fibre can be identified with the set of all $\epsilon_{12}\in(A\cap C)^*$ whose restriction to the subspace $(A\cap C\cap D)+(A\cap B\cap C)$ is prescribed. Hence the number of elements in this fibre is

$$|A\cap C|/(|A\cap C\cap D||A\cap B\cap C|/|A\cap B\cap C\cap D|).$$

We see that

$$((f_1\star f_2)\star f_3)((x,u),\epsilon)$$

$$=\sum_{y\in X_B,z\in X_C;(\epsilon_1,\epsilon_2,\epsilon_3)\in\mathfrak{T}'}f_1((x,y),\epsilon_1)f_2((y,z)\epsilon_2)f_3((z,u),\epsilon_3)\times$$

$$\frac{|A\cap C\cap D}{|A\cap D|}\frac{|A\cap B\cap C|}{|A\cap C|}\frac{|A\cap C|}{|A\cap C\cap D|}\frac{|A\cap B\cap C|}{|A\cap B\cap C\cap D|}$$

$$=\sum_{y\in X_B,z\in X_C;\epsilon_1,\epsilon_2,\epsilon_3}f_1((x,y),\epsilon_1)f_2((y,z),\epsilon_2)f_3((z,u),\epsilon_3)|A\cap B\cap C\cap D|/|A\cap D|.$$

A similar calculation applied to $(f_1\star(f_2\star f_3))((x,u),\epsilon)$ gives the same result. This proves (a).

From (a) we see that $\mathcal{F}\times\mathcal{F}\to\mathcal{F}$ makes $\mathcal{F}$ into an associative algebra over $\mathbf{C}$.

**1.3.** Let $A, B$ be in $\mathcal{L}$. Given an $E$-stable subset $\Xi$ of $X^2_{AB}$ and $\epsilon \in (A \cap B)^*$ we define $[\Xi]^\epsilon \in \mathcal{F}_{AB}$ to be the function which equals 1 on $((x,y),\epsilon)$ with $(x,y) \in \Xi$ and is equal to 0 at all other points of $X_{AB} \times (A \cap B)^*$. Note that

(a) *$\{[\mathcal{O}]^\epsilon; \mathcal{O} \in \overline{X^2}_{AB}, \epsilon \in (A \cap B)^*\}$ is a $\mathbf{C}$-basis $\underline{\mathcal{B}}_{AB}$ of $\mathcal{F}_{AB}$. Hence $\underline{\mathcal{B}} := \sqcup_{A,B}\underline{\mathcal{B}}_{AB}$ is a $\mathbf{C}$-basis of $\mathcal{F}$.*

Let $\Delta_A$ be the diagonal in $X^2_{AA}$. From the definitions we see that $\sum_{A \in \mathcal{L}}[\Delta_A]^0$ is a unit element of the algebra $\mathcal{F}$.

We have an isomorphism $\underline{K}_E(X^2) \to \mathcal{F}$ which maps the basis element corresponding to an irreducible $E$-equivariant vector bundle $V$ on $X^2$ to $[\mathcal{O}]^\epsilon$ where $\mathcal{O}$ is the set of all $(x,y) \in X^2$ such that the fibre of $V$ at $(x,y)$ is nonzero (hence one dimensional) and $\epsilon$ is the character by which the stabilizer of $(x,y)$ in $E$ acts on that fibre. This isomorphism is compatible with the algebra structures. (This provides another proof for the associativity of multiplication in $\mathcal{F}$.)

**1.4.** Let $B, C, D$ be in $\mathcal{L}$. Now $E \times C$ acts on $X^2_{BD}$ by

$$(e,c) : (x,y) \mapsto (e+x, e+c+y).$$

Let $\overline{X^2}_{BCD}$ be the set of $E \times C$-orbits on $X^2_{BD}$. The stabilizer of $(x,y) \in X^2_{BD}$ for the $E \times C$-action is $\{(e,c) \in E \times C; e \in B, e+b \in D\}$; this is a group isomorphic to

$$U_{BCD} := \{(b,c,d) \in B \times C \times D; b+c+d = 0\}$$

via $(e,c) \mapsto (e,c,e+c)$. Hence

(a) *if $\mathcal{M} \in \overline{X^2}_{BCD}$, then $|\mathcal{M}| = |E||C|/|U_{BCD}|$.*

**1.5.** Let $B, C, D$ be in $\mathcal{L}$. For any $\alpha \in (B \cap C \cap D)^*$ let $\tilde{\alpha}$ be the set of all $\epsilon \in (B \cap D)^*$ such that the restriction of $\epsilon$ to $B \cap C \cap D$ is equal to $\alpha$. Given an $E \times C$-stable subset $\mathcal{Y}$ of $X_{BD}$ and $\alpha \in (B \cap C \cap D)^*$ we define $[\mathcal{Y}]^\alpha = \sum_{\epsilon \in \tilde{\alpha}}[\mathcal{Y}]^\epsilon \in \mathcal{F}_{BD}$.

Let $\mathcal{B}_{BCD}$ be the collection consisting of the elements $[\mathcal{M}]^\alpha \in \mathcal{F}_{BD}$ for various $\mathcal{M} \in \overline{X^2}_{BCD}$ and various $\alpha \in (B \cap C \cap D)^*$. Let $\mathcal{B}_{CD} = \cup_{B \in \mathcal{L}}\mathcal{B}_{BCD}$ (this is a disjoint union). Clearly, $\mathcal{B}_{CD}$ is a linearly independent subset of $\mathcal{F}$. It spans a $\mathbf{C}$-subspace of $\mathcal{F}$ denoted by $[[CD]]$.

**Theorem 1.6.** *Let $C \in \mathcal{L}, D \in \mathcal{L}_X$.*

*(a) $[[CD]]$ is a left ideal of $\mathcal{F}$. Hence it is an $\mathcal{F}$-module.*

*(b) This $\mathcal{F}$-module is positive.*

## 2. Proof of Theorem 1.6

**2.1.** Let $A, B, D$ be in $\mathcal{L}_X$, $C \in \mathcal{L}$ and let $o' \in \bar{X}_B$. Let $\mathcal{O}' \in \overline{X^2}_{AB}, \mathcal{M} \in \overline{X^2}_{BCD}$. Assume that $p_2\mathcal{O}' = p_1\mathcal{O} = o'$. Let $\mathcal{O}' \bullet \mathcal{M} \subset X^2_{AD}$ be the image of

$$\widetilde{\mathcal{O}' \bullet \mathcal{M}} := \{(x,y,z) \in X_A \times X_B \times X_D; (x,y) \in \mathcal{O}', (y,z) \in \mathcal{M}\}$$

under the projection $X_A\times X_B\times X_D \to X_A\times X_D$, $(x,y,z)\mapsto(x,z)$. We have $|\mathcal{O}'| = |E|/|A\cap B|, |\mathcal{M}| = |E||C|/|U_{BCD}|$, $|o'| = |E|/|B|$. Since the second projection $\mathcal{O}' \to o'$ is surjective and $E$-equivariant, its fibre $F'_y$ at $y\in o'$ satisfies

$$|F'_y| = |\mathcal{O}'|/|o'| = |B|/|A\cap B|;$$

similarly, the fibre $F_y$ at $y\in o'$ of the first projection $\mathcal{M}\to o'$ satisfies

$$|F_y| = |\mathcal{M}|/|o'| = |B||C|/|U_{BCD}|.$$

It follows that

$$|\widetilde{\mathcal{O}'\bullet\mathcal{M}}| = \sum_{y\in o'}|F_y||F'_y| = \frac{|E|}{|B|}\frac{|B|}{|A\cap B|}\frac{|B||C|}{|U_{BCD}|} = \frac{|E||B||C|}{|A\cap B||U_{BCD}|}.$$

Now let $(x,z)\in\mathcal{O}'\bullet\mathcal{M}$. Let

$$Z := Z_{x,z} = \{y\in o'; (x,y)\in\mathcal{O}', (y,z)\in\mathcal{M}\},$$

(a nonempty set). Let $U = U_{A,B,C,D} = \{a,b,c,d)\in A\times B\times C\times D; a+b+c+d=0\}$. Now $U$ acts on $Z$ by $(a,b,c,d) : y\mapsto a+y = c+d+y$. Indeed, $(x,a+y) = (a+x,a+y)\in\mathcal{O}'$, $(c+d+y,z) = (c+d+y,d+z)\in\mathcal{M}$. This action is transitive: if $(x,y)\in\mathcal{O}'$, $(y,z)\in\mathcal{M}$, $(x,y')\in\mathcal{O}'$, $(y',z)\in\mathcal{M}$, $y\in o', y'\in o'$, then $y' = a+y$ for some $a\in A$ and $y' = c+d+y$ for some $c\in C$, $d\in D$. Then $b = a+c+d\in B$ so that $(a,b,c,d)\in U$. The stabilizer of $y\in o'$ in this action is

$$U' = \{(a,b,c,d)\in U; a+y = c+d+y = y\} = \{(a,b,c,d)\in U; a+y = a+b+y = y\}$$

that is

$$U' = \{(a,b,c,d)\in U; a\in A\cap B\}.$$

Hence we have a free transitive action of $U/U'$ on $Z$. Thus,

$$|Z_{x,z}| = |U|/|U'| = |U|/(|A\cap B||U_{BCD}|).$$

(We have an isomorphism $U' \xrightarrow{\sim} (A\cap B)\times U_{BCD}$ given by

$$(a,b,c,d)\mapsto(a,(a+b,c,d)).)$$

We have an injective homomorphism

$$(A\cap D)/(A\cap B\cap D)\to U/U', a = d\mapsto(a,0,0,d).$$

It follows that

$$\text{(a)}\qquad N_{A,B,C,D} := \frac{|U/U'|}{\frac{|A\cap D|}{|(A\cap B\cap D|}} \in \{1,2,2^2,2^3,\dots\}.$$

Since $|Z_{x,z}|$ is independent of $(x,z)$, it follows that

$$|\mathcal{O}' \bullet \mathcal{M}| = |\widetilde{\mathcal{O}' \bullet \mathcal{M}}|/(|U|/|U'|) = |E||B||C|/|U|,$$

so that

$$|\mathcal{O}' \bullet \mathcal{M}| = (|E||C|/|U_{ACD}|) \times (|B||U_{ACD}|/|U|)$$

We have a surjective map $U \to B \cap (A+C+D)$ given by $(a,b,c,d) \mapsto b$ whose kernel can be identified with $U_{ACD}$. Hence

$$|B||U_{ACD}|/|U| = |B|/|B \cap (A+C+D)| = |(A+B+C+D)/(A+C+D)|.$$

We see that

(b) $$|\mathcal{O}' \bullet \mathcal{M}| = (|E||C|/|U_{ACD}|) \times |(A+B+C+D)/(A+C+D)|.$$

Next we note that $\mathcal{O}' \bullet \mathcal{M}$ is a $E \times C$-stable subset of $X_{AD}$. (Indeed if $(x,z)$ is such that $(x,y) \in \mathcal{O}', (y,z) \in \mathcal{M}$ and $(e,c) \in E \times C$ then $(e+x, e+y) \in \mathcal{O}'$ and $(e+y, e+c+z) \in \mathcal{M}$.) Using (b) and 1.4(a) we see that

(c) $\mathcal{O}' \bullet \mathcal{M}$ is the disjoint union of $|(A+B+C+D)/(A+C+D)|$ orbits of $E \times C$ on $X_{AD}$.

**2.2.** We now assume that $\epsilon \in (A \cap B)^*$ and $\alpha \in (B \cap C \cap D)^*$. Let $\tilde{Y}$ be the set of all $(\epsilon', \epsilon_1) \in (A \cap D)^* \times (B \cap D)^*$ such that the restriction of $\epsilon_1$ to $B \cap C \cap D$ is $\alpha$ and the restrictions of $\epsilon', \epsilon, \epsilon_1$ to $A \cap B \cap D$ have sum 0.

Let $Y \subset (A \cap D)^*$ be the image of $\tilde{Y}$ under $(\epsilon', \epsilon_1) \mapsto \epsilon'$. From the definitions we have

(a) $$[\mathcal{O}']^\epsilon \star [\mathcal{M}]^\alpha = N_{A,B,C,D} \sum_{\mathcal{M}' \in \overline{X^2}_{ACD}; \mathcal{M}' \subset \mathcal{O}' \bullet \mathcal{M}} \sum_{\epsilon' \in Y} [\mathcal{M}']^{\epsilon'}.$$

We show:

(b) If $\epsilon' \in Y$ and $\delta \in (A \cap D)^*$ has restriction 0 to $A \cap C \cap D$ then $\epsilon' + \delta \in Y$.

Let $\epsilon_1 \in (B \cap D)^*$ be an element whose restriction to $B \cap C \cap D$ is equal to $\alpha$ and such that the restrictions of $\epsilon', \epsilon, \epsilon_1$ to $A \cap B \cap D$ have sum 0.

It is enough to show that there exists $\tilde{\epsilon} \in (B \cap D)^*$ whose restriction to $B \cap C \cap D$ is 0 and such that the restrictions of $\epsilon' + \delta, \epsilon, \epsilon_1 + \tilde{\epsilon}$ to $A \cap B \cap D$ have sum 0 (or equivalently, the restrictions of $\delta, \tilde{\epsilon}$ to $A \cap B \cap D$ are equal). Let $\delta' : A \cap B \cap D \to \mathbf{Z}/2$ be the restriction of $\delta$. We set $B' = B \cap D, A' = A \cap B \cap D, C' = B \cap C \cap D$. We have $A' + C' \subset B'$ We must prove that there exists $\tilde{\epsilon} : B' \to \mathbf{Z}/2$ (linear) whose restriction to $C'$ is zero, such that the restriction of $\tilde{\epsilon}$ to $A'$ equals $\delta'$. We can find a subspace $A''$ of $A'$ and a subspace $B''$ of $B'$ such that $A' = (A' \cap C') \oplus A''$, $B' = (A' + C') \oplus B''$. We have $A' + C' = A'' \oplus C'$. Note that the restriction of $\delta'$ to $A' \cap C'$ is 0. We define $\tilde{\epsilon} : B' \to \mathbf{Z}/2$ (linear) to be $\delta'$ on $A''$ and to be zero on $B''$ and $C'$. This satisfies our requirement. This proves (b).

Let $Y'$ be the image of $Y$ under the obvious surjective map $(A\cap D)^*\to(A\cap C\cap D)^*$. From (b) we see that $Y$ is a union of fibres of this last map. Hence (a) can be rewritten as follows.

(c) $$[\mathcal{O}']^\epsilon\star[\mathcal{M}]^\alpha=N_{A,B,C,D}\sum_{\mathcal{M}'\in\overline{X^2}_{ACD};\mathcal{M}'\subset\mathcal{O}'\bullet\mathcal{M}}\sum_{\alpha'\in Y'}[\mathcal{M}']^{\alpha'}.$$

From this we see that 1.6(a),(b) hold.

**2.3.** Let $C\in\mathcal{L},D\in\mathcal{L}_X$ We have a direct sum decomposition

$$[[CD]]=\oplus_{o\in\bar{X}_D}[[Co]]$$

where $[[Co]]$ is the subspace of $[[CD]]$ spanned by $\mathcal{B}_{Co}=\{[\mathcal{M}]^\alpha\in\mathcal{B}_{CD};p_2\mathcal{M}=o\}$. From 1.6 we see that for any $o\in\bar{X}_D$:

(a) $[[Co]]$ is a left ideal of $\mathcal{F}$ and the $\mathcal{F}$-module $[[Co]]$ is positive.

We have the following result.

**Theorem 2.4.** *Let $C,D$ be in $\mathcal{L}_X$. Let $\mathcal{O}\in\overline{X^2}_{CD}$ be such that $p_2\mathcal{O}=o$ and let $\alpha\in(C\cap D)^*$. We have $[\mathcal{O}]^\alpha\in[[Co]]$. Moreover $[\mathcal{O}]^\alpha$ is a generator of the $\mathcal{F}$-module $[[Co]]$.*

We have also $\mathcal{O}\in\overline{X^2}_{CCD}$ hence 2.2(c) with $\mathcal{O}$ replacing $\mathcal{M}$ is applicable.

In this case we have by 2.1(a): $|\mathcal{O}'\bullet\mathcal{O}|=|E||C|/|U_{ACD}|$ since $A+C+C+D=A+C+D$. Using 1.4(a) we see that any $\mathcal{M}'$ in 2.2(c) satisfies $|\mathcal{M}'|=|E||C|/|U_{ACD}|$ that is $|\mathcal{M}'|=|\mathcal{O}'\bullet\mathcal{O}|$. It follows that

(a) $\mathcal{O}'\bullet\mathcal{O}$ is a single $E\times C$-orbit.

We have $\epsilon\in(A\cap C)^*$, $\alpha\in(C\cap D)^*$ and $Y\subset(A\cap D)^*$ is the set of all $\epsilon'\in(A\cap D)^*$ such that the restrictions of $\epsilon',\epsilon,\alpha$ to $A\cap C\cap D$ have sum 0. Now $Y'$ consists of a single element $\alpha'$. We see that

$$[\mathcal{O}']^\epsilon\star[\mathcal{O}]^\alpha=N_{A,C,C,D}\mathcal{M}'^{\alpha'}$$

where $\mathcal{M}'=\mathcal{O}'\bullet\mathcal{O}\in\in\overline{X^2}_{ACD},p_2\mathcal{M}'=o$ and $\alpha'$ is the restriction of some $\epsilon'$ as above to $A\cap C\cap D$.

Conversely, let $\mathcal{M}'\in\overline{X^2}_{ACD}$ such that $p_2\mathcal{M}'=o$ and let $\alpha'$ be a fibre of the obvious map $(A\cap D)^*\to(A\cap C\cap D)^*$. Let $(x,z)\in\mathcal{M}'$. Since $z\in o=p_2\mathcal{O}$ we can find $y\in X_C$ such that $(y,z)\in\mathcal{O}$. Let $\mathcal{O}'\in\overline{X^2}_{AC}$ be the $E$-orbit of $(x,y)$. Then $\mathcal{O}'\bullet\mathcal{O}$ is defined and is an $E\times C$-orbit. Since $\mathcal{O}'\bullet\mathcal{O},\mathcal{M}'$ are $E\times C$-orbits and both contain $(x,z)$, we must have $\mathcal{M}'=\mathcal{O}'\bullet\mathcal{O}$. Let $\epsilon'\in(A\cap D)^*$ be such that the restriction of $\epsilon'$ to $A\cap C\cap D$ is $\alpha'$. We can find $\alpha\in(C\cap D)^*$ such that the restrictions of $\epsilon',\epsilon,\alpha$ to $A\cap C\cap D$ have sum 0. We have

$$[\mathcal{O}']^\epsilon\star[\mathcal{O}]^\alpha=N_{A,C,C,D}\mathcal{M}'^{\alpha'}.$$

We now see that $[[Co]]$ is exactly the left ideal of $\mathcal{F}$ generated by the element $[\mathcal{O}]^\alpha$. This proves the theorem.

## 3. Trace computation

**3.1.** Let $C\in\mathcal{L},D\in\mathcal{L}_X$ and let $o\in\bar{X}_D$. Let $A\in\mathcal{L}_X,B\in\mathcal{L}_X$. The following result gives the trace $t$ of an element $[\mathcal{O}']^\epsilon\in\mathcal{F}$ (with $\mathcal{O}'\in\overline{X^2}_{AB},\epsilon\in(A\cap B)^*$) acting in the $\mathcal{F}$-module $[[Co]]$.

**Theorem 3.2.** *(a) If $p_1\mathcal{O}'\ne p_2\mathcal{O}'$ then $t=0$.*

*(b) Assume now that $p_1\mathcal{O}'=p_2\mathcal{O}'=o'$ so that $A=B$. We can find $e\in E$ such that $\mathcal{O}'=\{(x,e+x);x\in o'\}$. Let $\delta=\delta_{e\in B+C+D}$, $\delta'=\delta_{\epsilon|_{B\cap C\cap D}=0}$. We have $t=\delta\delta'|B\cap C\cap D||E|/|B+C+D|$.*

We will compute $t$ using the matrix by which $[\mathcal{O}']^\epsilon$ acts in terms of the basis $\mathcal{B}_{Co}$. Then it is clear that (a) holds. In the rest of the proof we assume that we are in the setup of (b).

Let $\mathcal{M}\in\overline{X^2}_{BCD}$ with $p_2\mathcal{M}=o$, $p_1\mathcal{M}=o'$ and let $\alpha\in(B\cap C\cap D)^*$. We show (denoting by $\bar\epsilon$ the restriction of $\epsilon$ to $B\cap C\cap D$):

(c) $[\mathcal{O}']^\epsilon\star[\mathcal{M}]^\alpha=[\mathcal{O}'\bullet\mathcal{M}]^{\alpha+\bar\epsilon}$ with $\mathcal{O}'\bullet\mathcal{M}\in\overline{X^2}_{BCD}$.

We will use 2.1(c) with $A=B$. In our case from 2.1(a) we see that $N_{A,B,C,D}=1$. Note that in 2.1(c) we have $|(A+B+C+D)/(A+C+D)|=1$ so that in our case $\mathcal{O}'\bullet\mathcal{M}$ is a single $E\times C$-orbit on $X_{AD}=X_{BD}$. We have $\epsilon\in(A\cap B)^*=B^*$ and $\alpha\in(B\cap C\cap D)^*=(A\cap C\cap D)^*$. In our case the set $Y$ of 2.2 is the set of all $\epsilon'\in(B\cap D)^*$ such that the restrictions $\bar\epsilon,\bar\epsilon'$ of $\epsilon,\epsilon'$ to $B\cap C\cap D$ satisfy $\bar\epsilon+\bar\epsilon'+\alpha=0$. Hence the set $Y'$ of 2.2 consists of a single element, namely $\alpha+\bar\epsilon$. We see that (c) holds.

We set ${}^e\mathcal{M}=\{(x,z);(e+x,z)\in\mathcal{M}\}\subset X^2_{BD}$. It is clear that ${}^e\mathcal{M}\in\overline{X^2}_{BCD},p_2{}^e\mathcal{M}=o$, and that $\mathcal{O}'\bullet\mathcal{M}={}^e\mathcal{M}$. From the definitions we see that the condition that ${}^e\mathcal{M}=\mathcal{M}$ is that $e\in B+C+D$. Moreover the condition that $\alpha=\alpha+\bar\epsilon$ is that $\epsilon|_{B\cap C\cap D}=0$. We now see that

$$t=\delta\delta'|\{\mathcal{M}\in\overline{X^2}_{BCD};p_2\mathcal{M}=o,p_1\mathcal{M}=o'\}||(B\cap C\cap D)^*|,$$

$$\begin{aligned}t&=\delta\delta'|\overline{X^2}_{BCD}||B\cap C\cap D|/(|\bar{X}_B||\bar{X}_D|)\\&=\delta\delta'|B\cap C\cap D||X_B||X_D|/(|E||C|/|U_{BCD}||\bar{X}_B||\bar{X}_D|)\\&=\delta\delta'|B\cap C\cap D|(|E|/|B|)(|E|/|D|)/(|E||C|/|U_{BCD}|)\\&=\delta\delta'|B\cap C\cap D|(|E||U_{BCD}|/(|B||C||D|).\end{aligned}$$

We have $|U_{BCD}|=|(C+D)\cap B||C\cap D|$ hence

$$\begin{aligned}t&=\delta\delta'|B\cap C\cap D||E||(C+D)\cap B||C\cap D|/(|B||C+D||C\cap D|)\\&==\delta\delta'|B\cap C\cap D||E||(C+D)\cap B|/(|B||C+D|)\\&=\delta\delta'|B\cap C\cap D||E|/|B+C+D|.\end{aligned}$$

This completes the proof.

**Corollary 3.3.** *Let $C\in\mathcal{L},C'\in\mathcal{L},D\in\mathcal{L}_X,D'\in\mathcal{L}_X$ and let $o\in\bar{X}_D$, $o'\in\bar{X}_{D'}$. Assume that $C\cap D=C'\cap D'$, $C+D=C'+D'$. Then for any $f\in\mathcal{F}$ we have $\mathrm{tr}(f,[[Co]])=\mathrm{tr}(f',[[C'o']])$.*

## 4. The algebra $\bar{\mathcal{F}}^e$

**4.1.** Let $e \in E$. Let $X^e = \{x \in X; e + x = x\}$, $(X^2)^e = X^e \times X^e$. Then $X^e$ is $E$-stable and $(X^2)^e \subset X^2$ is $E$-stable. Let $\mathcal{L}^e = \{A \in \mathcal{L}; e \in A\}$, $\mathcal{L}^e_X = \mathcal{L}_X \cap \mathcal{L}^e$. We have $X^e = \sqcup_{A \in \mathcal{L}^e} X_A$ and $(X^2)^e = \sqcup_{A,B \text{ in } \mathcal{L}^e} X^2_{AB}$.

Let $\bar{\mathcal{F}}^e$ be the set of functions $(X^2)^e \to \mathbf{C}$ that are constant on each $E$-orbit. We have $\bar{\mathcal{F}}^e = \oplus_{A,B \text{ in } \mathcal{L}} \bar{\mathcal{F}}^e_{AB}$ where, if $e \in A \cap B$, $\bar{\mathcal{F}}^e_{AB}$ is the subspace of $\bar{\mathcal{F}}^e$ consisting of functions with support contained in $X^2_{AB}$, while if $e \notin A \cap B$, we have $\bar{\mathcal{F}}^e_{AB} = 0$.

We define a $\mathbf{C}$-bilinear map $\bar{\mathcal{F}}^e \times \bar{\mathcal{F}}^e \to \bar{\mathcal{F}}^e$, $(f_1, f_2) \mapsto f_1 \sharp f_2$ by

$$(f_1 \sharp f_2)(x, z) = \sum_{y \in X^e} f_1(x, y) f_2(y, z)$$

for $(x, z) \in (X^2)^e$. This defines an associative algebra structure on $\bar{\mathcal{F}}^e$ with unit element being the function which is 1 on any $(x, x)$ with $x \in X^e$ and is zero on all other points of $(X^2)^e$. This algebra can be identified with the algebra of $E$-equivariant endomorphisms of the permutation $E$-module of functions $X^e \to \mathbf{C}$; hence it is a semisimple algebra.

For $A, B$ in $\mathcal{L}$ we define a $\mathbf{C}$-linear map

$$\tau^e = \tau^e_{AB} : \mathcal{F}_{AB} \to \bar{\mathcal{F}}^e_{AB}$$

to be 0 if either $e \notin A$ or $e \notin B$ and to be $f \mapsto \tau^e(f)$ with

$$\tau^e(f)(x, y) = \sum_{\epsilon \in (A \cap B)^*} f((x, y), \epsilon)(-1)^{\epsilon(e)}$$

if $e \in A \cap B$. Taking direct sum over all $A, B$ we obtain a $\mathbf{C}$-linear map $\mathcal{F} \to \bar{\mathcal{F}}^e$, $f \mapsto \tau^e(f)$.

For $A, B, C$ in $\mathcal{L}$ and for $f_1 \in \mathcal{F}_{AB}$, $f_2 \in \mathcal{F}_{BC}$ we show:

(a) $$\tau^e(f_1 \star f_2) = \tau^e(f_1) \sharp \tau^e(f_2).$$

Assume first that $e \in A \cap C$. Let $(x, z) \in (X^2)^e_{AC}$. We have

$$\tau^e(f_1 \star f_2)(x, z)$$
$$= \sum_{(y, \epsilon_1, \epsilon_2, \epsilon) \in X_B \times (A \cap B)^* \times (B \cap C)^* \times (A \cap C)^*; \bar{\epsilon}_1 + \bar{\epsilon}_2 + \bar{\epsilon} = 0}$$
$$(-1)^{\epsilon(e)} f_1((x, y), \epsilon_1) f_2((y, z), \epsilon_2) |A \cap B \cap C| / |A \cap C|$$
$$= \sum_{(y, \epsilon_1, \epsilon_2) \in X_B \times (A \cap B)^* \times (B \cap C)^*}$$
$$f_1((x, y), \epsilon_1) f_2((y, z), \epsilon_2) (|A \cap B \cap C| / |A \cap C|) z(\epsilon_1, \epsilon_2)$$

where

$$z_{\epsilon_1,\epsilon_2} \sum_{\epsilon\in(A\cap C)^*;\bar{\epsilon}=\bar{\epsilon}_1+\bar{\epsilon}_2} (-1)^{\epsilon(e)}.$$

We now assume further that $e \in B$, so that $e \in A \cap B \cap C$; then for each $\epsilon$ in the last sum we have $\epsilon(e) = \bar{\epsilon}(e) = \bar{\epsilon}_1(e) + \bar{\epsilon}_2(e) = \epsilon_1(e) + \epsilon_2(e)$ so that

$$z_{\epsilon_1,\epsilon_2} = \sum_{\epsilon\in(A\cap C)^*;\bar{\epsilon}=\bar{\epsilon}_1+\bar{\epsilon}_2} (-1)^{\epsilon_1(e)+\epsilon_2(e)}(|A\cap C|/|A\cap B\cap C|)(-1)^{\epsilon_1(e)}(-1)^{\epsilon_2(e)}.$$

We see that

$$\begin{aligned}
&\tau^e(f_1\star f_2)(x,z)\\
&=\sum_{y\in X_B}\sum_{\epsilon_1\in(A\cap B)^*} f_1((x,y),\epsilon_1)(-1)^{\epsilon_1(e)} \sum_{\epsilon_2\in(B\cap C)^*} f_2((y,z),\epsilon_2)(-1)^{\epsilon_2(e)}\\
&=\sum_{y\in X_B}\tau^e(f_1)(x,y)\tau^e(f_2)(y,z)=(\tau^e(f_1)\sharp\tau^e(f_2)(x,z),
\end{aligned}$$

as desired.

Next we assume that $e \in A\cap C$, $e \notin B$. We show that $z(\epsilon_1,\epsilon_2) = 0$ for any $e_1,\epsilon_2$. It is enough to show that

(b) the homomorphism $\epsilon \mapsto \epsilon(e)$ from the group

$$\{\epsilon\in(A\cap C)^*;\bar{\epsilon}=0\}=\{\epsilon\in(A\cap C)^*;\epsilon=0\text{ on }A\cap B\cap C\}$$

to $\mathbf{Z}/2$ is nontrivial.

Let $V = A\cap C, V' = V\cap B \subset V$. We have $e \in V - V'$. We write $V = V' \oplus Z$ where $Z$ is a subspace containing $e$. We can find $\epsilon \in V^*$ such that $\epsilon(e) \neq 0$, $\epsilon|V' = 0$. Thus (b) holds. We see that $\tau^e(f_1\star f_2)(x,z) = 0$. On the other hand we have $\tau^e(f_1) = 0, \tau^e(f_2) = 0$. Thus (a) holds in this case.

If $e \notin A$ or $e \notin C$ then both sides of (a) are zero. This proves (a).

**4.2.** For $A, B$ in $\mathcal{L}$ we show:

(a) the $\mathbf{C}$-linear map $\tau = \tau_{AB} : \mathcal{F}_{AB} \to \oplus_{e\in E}\bar{\mathcal{F}}^e_{AB}$ given by $f \mapsto (\tau^e(f); e\in E)$ is an isomorphism.

Let $f \in \mathcal{F}_{AB}$ be in the kernel of $\tau$. Then for any $(x,y) \in X^2_{AB}$ we have $f((x,y),\epsilon)(-1)^{\epsilon(e)} = 0$ for any $e \in A\cap B$. Then the function $\epsilon \mapsto f((x,y),\epsilon)$ on $(A\cap B)^*$ has Fourier transform (a function on $A\cap B$) equal to 0. Since the Fourier transform is an isomorphism, it follows that $f = 0$. This proves the injectivity of the map (a). It remains to show that (a) is a linear map betweem two vector spaces of the same dimension; but these dimensions are $|\overline{X^2}_{AB}||(A\cap B)^*|$, $|\overline{X^2}_{AB}||A\cap B|$ and these are obviously equal. This proves (a).

Taking direct sum over all $A, B$ in $\mathcal{L}$ we obtain an isomorphism

(b) $\tau : \mathcal{F} \to \oplus_{e\in E}\bar{\mathcal{F}}^e$

which by 4.1(a) is an algebra isomorphism. It follows that

(c) $\mathcal{F}$ is a semisimple algebra.

**4.3.** Let $e,e'$ in $E$. We define a linear map $T_{e'}:\bar{\mathcal{F}}^e\to\bar{\mathcal{F}}^e$, $f\mapsto T_{e'}f$ by $(T_{e'}f)(x,y)=f(e'+x,y)=f(x,e'+y)$. Now $e'\mapsto T_{e'}$ is a homomorphism $E\to GL(\bar{\mathcal{F}}^e)$. Hence we have a direct sum decomposition $\bar{\mathcal{F}}^e=\oplus_{e^*\in E^*}\bar{\mathcal{F}}^{e,e^*}$ where

$$\bar{\mathcal{F}}^{e,e^*}=\{f\in\bar{\mathcal{F}}^e;T_{e'}f=(-1)^{e^*(e')}f\quad\forall e'\in E\}.$$

Let $e^*\in E^*$, $f_1\in\bar{\mathcal{F}}^e$, $f_2\in\bar{\mathcal{F}}^{e,e^*}$. We show that

(a) $$f_1\sharp f_2\in\bar{\mathcal{F}}^{e,e^*}$$

. Indeed, for $e'\in E$ and $(x,z)\in(X^2)^e$ we have

$$\begin{aligned}&T_{e'}(f_1\sharp f_2)(x,z)=(f_1\sharp f_2)(x,e'+z)\\&=\sum_{y\in X^e}f_1(x,y)f_2(y,e'+z)=\sum_{y\in X^e}(-1)^{e^*(e')}f_1(x,y)f_2(y,z)\\&=(-1)^{e^*(e')}(f_1\sharp f_2)(x,z).\end{aligned}$$

This proves (a).

We now show that

(b) $$f_2\sharp f_1\in\bar{\mathcal{F}}^{e,e^*}$$

. Indeed, for $e'\in E$ and $(x,z)\in(X^2)^e$ we have

$$\begin{aligned}&T_{e'}(f_2\sharp f_1)(x,z)=(f_2\sharp f_1)(e'+x,z)\\&=\sum_{y\in X^e}f_1(e'+x,y)f_2(y,z)=\sum_{y\in X^e}(-1)^{e^*(e')}f_1(x,y)f_2(y,z)\\&=(-1)^{e^*(e')}(f_1\sharp f_2)(x,z).\end{aligned}$$

This proves (b).

From (a),(b) we see that

(c) for any $e^*\in E^*$, $\bar{\mathcal{F}}^{e,e^*}$ is a two-sided ideal of $\bar{\mathcal{F}}^e$.

We show that

(d) if $e^*\in E^*$, $\bar{\mathcal{F}}^{e,e^*}\neq0$ then $e^*(e)=0$.

Indeed, let $f\in\bar{\mathcal{F}}^{e,e^*}-\{0\}$ and let $(x,y)\in(X^2)^e$ be such that $f(x,y)\neq0$. We have $f(x,y)=f(e+x,y)=(-1)^{e^*(e)}f(x,y)$. Since $f(x,y)\neq0$, it follows that $(-1)^{e^*(e)}=1$ so that $e^*(e)=0$, as desired.

For $(e,e^*)\in E\times E^*$ we set

(e) $$\mathcal{L}_X^{e,e^*}=\{A\in\mathcal{L}_X^e;e^*|_A=0\}.$$

Note that if $\mathcal{L}_X^{e,e^*}\neq\emptyset$ then $e^*(e)=0$.

**4.4.** Let $e \in E$. Let $\overline{(X^2)^e}$ be the set of $E$-orbits in $(X^2)^e$. For any $\mathcal{O} \in \overline{(X^2)^e}$ we denote by $[\mathcal{O}]$ the function $(X^2)^e \to \mathbf{C}$ which equals 1 on $\mathcal{O}$ and equals 0 on $(X^2)^e - \mathcal{O}$. Then $\{[\mathcal{O}]; \mathcal{O} \in \overline{(X^2)^e}\}$ is a basis of $\bar{\mathcal{F}}^e$.

**Proposition 4.5.** *Let $e \in E, e^* \in E^*$ be such that $e^*(e) = 0$. Let $A, B$ be in $\mathcal{L}_X^e$. For $\mathcal{O} \in \overline{(X^2)^e}$ with $\mathcal{O} \subset X_A \times X_B$, let $t$ be the trace of left multiplication by $[\mathcal{O}]$ on $\bar{\mathcal{F}}^{e,e^*}$. If $p_1\mathcal{O} \ne p_2\mathcal{O}$ or if $e^*|_A \ne 0$, then $t = 0$. If $e^*|_A = 0$ and $p_1\mathcal{O} = p_2\mathcal{O} = o$ (so that $A = B$ and for some $e_0 \in E$ we have $\mathcal{O} = \{(e_0 + x, x), x \in o\}$), then*

$$t = (-1)^{e^*(e_0)} \sum_{D \in \mathcal{L}_X^{e,e^*}} |\bar{X}_D|.$$

Note that the projection $\pi^{e,e^*}$ of $\bar{\mathcal{F}}^e$ onto its direct summand $\bar{\mathcal{F}}^{e,e^*}$ is an algebra homomorphism; it is given by

$$f \mapsto \pi^{e,e^*}(f) = |E|^{-1} \sum_{e' \in E} (-1)^{e^*(e')} T_{e'}(f).$$

Let $t'$ be the trace of left multiplication by $\pi^{e,e^*}([\mathcal{O}])$ on $\bar{\mathcal{F}}^e$. We have

$$t = t' = \sum_{\mathcal{O}' \in \overline{(X^2)^e}} n_{\mathcal{O}',\mathcal{O}'}$$

where

$$|E|^{-1} \sum_{e' \in E} (-1)^{e^*(e')} T_{e'}([\mathcal{O}]) \sharp [\mathcal{O}'] = \sum_{\mathcal{O}'' \in \overline{(X^2)^e}} n_{\mathcal{O}'',\mathcal{O}'} [co''].$$

For $e' \in E$ we have $T_{e'}([\mathcal{O}]) = f_{e'[\mathcal{O}]}$ where

$${}^{e'}[\mathcal{O}] = \{(x, y); (e' + x, y) \in \mathcal{O}\} \in \overline{(X^2)^e}\}.$$

Hence for $\mathcal{O}' \in \overline{(X^2)^e}$ with $\mathcal{O}' \in X_C \times X_D$, $n_{\mathcal{O}',\mathcal{O}'}$ is the coefficient of $[\mathcal{O}']$ in

$$|E|^{-1} \sum_{e' \in E} (-1)^{e^*(e')} [{}^{e'}\mathcal{O}] \sharp [\mathcal{O}'].$$

Note that $n_{\mathcal{O}',\mathcal{O}'} = 0$ unless $p_2\mathcal{O} = p_1\mathcal{O} = p_1\mathcal{O}'$ (so that $A = B = C$) which we now assume. From the definitions we see that $[{}^{e'}\mathcal{O}] \sharp [\mathcal{O}'] = [{}^{e'+e_0}\mathcal{O}']$. The condition that $[{}^{e'+e_0}\mathcal{O}'] = [\mathcal{O}']$ is that $e' \in e_0 + A + D$. We see that

$$n_{\mathcal{O}',\mathcal{O}'} = |E|^{-1} \sum_{e' \in e_0 + A + D} (-1)^{e^*(e')}.$$

This is 0 if $e^*|_{A+D}\ne 0$ and is $|E|^{-1}(-1)^{e^*(e_0)}|A+D|$ if $e^*|_{A+D}=0$.

We see that if $p_1\mathcal{O}\ne p_2\mathcal{O}$ then $t=0$, while if $p_1\mathcal{O}=p_2\mathcal{O}$, then

$$t=|E|^{-1}(-1)^{e^*(e_0)}\sum_{D\in\mathcal{L}^e;e^*|_{A+D}=0}|A+D|$$
$$|\{\mathcal{O}'\in\overline{X^2}_{AD};p_1\mathcal{O}'=p_1\mathcal{O}\}|$$

that is

$$t=|E|^{-1}(-1)^{e^*(e_0)}\sum_{D\in\mathcal{L}^e;e^*|_{A+D}=0}|A+D||\bar{X}_D||E||A\cap D|/|A||D|$$
$$=(-1)^{e^*(e_0)}\sum_{D\in\mathcal{L}^e;e^*|_{A+D}=0}|\bar{X}_D|.$$

Thus $t=0$ if $A\notin\mathcal{L}_X^{e,e^*}$ and

$$t=(-1)^{e^*(e_0)}\sum_{D\in\mathcal{L}_X^{e,e^*}}|\bar{X}_D|$$

if $A\in\mathcal{L}_X^{e,e^*}$. The proposition is proved.

**Corollary 4.6.** *Let $e\in E,e^*\in E^*$ be such that $e^*(e)=0$. We have*

$$\dim\bar{\mathcal{F}}^{e,e^*}=(\sum_{A\in\mathcal{L}_X^{e,e^*}}|\bar{X}_A|)^2.$$

The unit element of $\bar{\mathcal{F}}^e$ is $\sum_{\mathcal{O}}[\mathcal{O}]$ where $\mathcal{O}$ is of the form $\tilde{o}=\{(x,x);x\in o\}$ for various $E$-orbits $o$ in $X^e$. Hence

$$\dim\bar{\mathcal{F}}^{e,e^*}=\sum_{\mathcal{O}}\operatorname{tr}([\mathcal{O}],\bar{\mathcal{F}}^{e,e^*})$$
$$=\sum_{A\in\mathcal{L}_X^{e,e^*}}|\{\tilde{o};o\subset X_A\}|\sum_{D\in\mathcal{L}_X^{e,e^*}}|\bar{X}_D|=(\sum_{A\in\mathcal{L}_X^{e,e^*}}|\bar{X}_A|)(\sum_{D\in\mathcal{L}_X^{e,e^*}}|\bar{X}_D|).$$

The corollary follows.

**4.7.** Let $M\subset P$ be subspaces of $E$ such that for some $D\in\mathcal{L}_X$ we have $M\subset D\subset P$ (we then have $M=C\cap D,P=C+D$ for some $C\in\mathcal{L}$.) Let $o\in\bar{X}_D$. We set $(MP)=[[Co]]$. From 3.3 we see that the trace of any element of $\mathcal{F}$ on the $\mathcal{F}$-module $(MP)$ is independent of the choice of $D,C,o$ above. Since $\mathcal{F}$ is a semisimple algebra, it follows that the isomorphism class of the $\mathcal{F}$-module $(MP)$ is independent of the choice of $D,C,o$ above. From 2.3(a) we see that

(a) The $\mathcal{F}$-module $(MP)$ is positive.

## 5. The centre of $\bar{\mathcal{F}}^e$

**5.1.** Let $e \in E$. Let $\bar{z}^e$ be the center of the algebra $\bar{\mathcal{F}}^e$. Let $\bar{z}'^e$ be the set of all $f \in \bar{\mathcal{F}}^e$ such that (a),(b),(c) below hold.

(a) The support of $f$ is contained in $\cup_{e' \in E}\{(x, e'+x); x \in X^e\}$

(b) For any $A \in \mathcal{L}_X^e$ there is a well defined function $f_A : E \to \mathbf{C}$ such that $f(x, e'+x) = f_A(e')$ for any $x \in X_A, e' \in E$. Note that $f_A(e'+a) = f_A(e')$ for $e' \in E, a \in A$.

(c) For any $A, B$ in $\mathcal{L}_X^e$ and any $e' \in E$ we have

$$\sum_{b \in B} f_A(e'+b) = \sum_{a \in A} f_B(e'+a).$$

**Proposition 5.2.** *We have $\bar{z}^e = \bar{z}'^e$.*

Let $f \in \bar{\mathcal{F}}^e$. We have $f \in \bar{z}^e$ if and only if for any $A, B$ in $\mathcal{L}_X^e$ and any $\mathcal{O} \in \overline{X^2}_{AB}$ we have

$$\sum_{y \in X^e} f(x,y) f_{\mathcal{O}}(y,z) = \sum_{y \in X^e} f_{\mathcal{O}}(x,y) f(y,z)$$

for any $(x,z) \in (X^2)^e$ or equivalently

$$\sum_{e_1 \in E_{A\cap B}} f(x, e_1+\alpha)\delta_{z=e_1+\beta} = \sum_{e'_1 \in E_{A\cap B}} \delta_{x=e'_1+\alpha} f(e'_1+\beta, z)$$

for any $x, z$ in $X^e$ and any $\alpha \in X_A, \beta \in X_B$ or equivalently

(a) $$\sum_{e_1 \in E} f(x, e_1+\alpha)\delta_{z=e_1+\beta} = \sum_{e'_1 \in E} \delta_{x=e'_1+\alpha} f(e'_1+\beta, z).$$

(For a subspace $C \subset E$ we let $E_C$ be a subspace of $E$ such that $E = C \oplus E_C$.)

Assume now that $f$ satisfies (a). If $\alpha = \beta = z$, $x \ne e'_1 + \alpha$ for any $e'_1 \in E$, then (a) becomes

$$\sum_{e_1 \in E} f(x, e_1+\alpha)\delta_{\alpha=e_1+\alpha} = 0.$$

Hence $\sum_{e_1 \in A} f(x, e_1+\alpha) = 0$ that is $|A| f(x,\alpha) = 0$ and $f(x,\alpha) = 0$. We see that $f$ satisfies 5.1(a).

If in (a) we have $x \ne e_1 + \alpha$ for any $\epsilon_1 \in E$, then the left hand side of (a) is 0 by 5.1(a) and the right hand side of (a) is also 0; if we have $z \ne e_2 + \beta$ for any $e_2 \in E$ the the left hand side of (a) is 0 and the right hand side of (a) is 0 by 5.1(a). Thus, in the presence of 5.1(a) the condition (a) is equivalent to:

$$\sum_{e_1 \in E} f(e_0+\alpha, e_1+\alpha)\delta_{e'_0+\beta=e_1+\beta} = \sum_{e'_1 \in E} \delta_{e_0+\alpha=e'_1+\alpha} f(e'_1+\beta, e'_0+\beta)$$

for any $\alpha\in X_A,\beta\in X_B$ and any $e_0,e'_0$ in $E$. Setting $b=e_1+e'_0\in B, a=e'_1+e_0\in A, e=e_0+e'_0$, this is equivalent to the condition that

(b) $$\sum_{b\in B} f(\alpha,e+b+\alpha)=\sum_{a\in A} f(\beta,e+a+\beta)$$

for any $\alpha\in X_A,\beta\in X_B$ and any $e\in E$. When $A=B$ this condition is

$$|B|f(\alpha,e+\alpha)=|A|f(\beta,e+\beta)$$

for any $\alpha,\beta$ in $X_A$ and any $e\in E$. Hence $f$ satisfies 5.1(b). Then (b) becomes 5.1(c). The proposition follows.

Let $(E\times E^*)_X$ be the set of all $(e,e^*)\in E\times E^*$ such that $\mathcal{L}_X^{e,e^*}\neq\emptyset$ (or equivalently that $\bar{\mathcal{F}}^{e,e^*}\neq 0$, see 4.6).

**Theorem 5.3.** *If $(e,e^*)\in(E\times E^*)_X$, then the algebra $\bar{\mathcal{F}}^{e,e^*}$ is simple.*

Since the algebra $\bar{\mathcal{F}}^e$ is semisimple, it is enough to show that

(a) $\dim(\bar{z}^e\cap\bar{\mathcal{F}}^{e,e^*})=1$.

Let $f\in\bar{z}^e\cap\bar{\mathcal{F}}^{e,e^*}$. For $e'\in E$ we have $T_{e'}f=(-1)^{e^*(e')}f$. Evaluating both sides at $(x,e_0+x)$ with $x\in X_A, A\in\mathcal{L}_X^e, e_0\in E$ we obtain $f(e'+x,e_0+x)=(-1)^{e^*(e')}f(x,e_0+x)$ that is $f_A(e_0+e')=(-1)^{e^*(e')}f_A(e_0)$ (with $f_A$ as in 5.1). In particular, $f_A(e')=(-1)^{e^*(e')}f_A(0)$. If now $A,B$ are in $\mathcal{L}_X^e$, then by 5.1(c),5.2, we have

$$\sum_{b\in B} f_A(b)=\sum_{a\in A} f_B(a),$$

that is,

$$\sum_{b\in B}(-1)^{e^*(b)}f_A(0)=\sum_{a\in A}(-1)^{e^*(a)}f_B(0)$$

or equivalently

(b) $$|B|\delta_{B\in\mathcal{L}_X^{e,e^*}}f_A(0)=|A|\delta_{A\in\mathcal{L}_X^{e,e^*}}f_B(0).$$

We see that if $A,B$ are in $\mathcal{L}_X^{e,e^*}$ then $f_A(0)/|A|=f_B(0)/|B|$. Thus for some $c\in\mathbf{C}$ we have $f_A(0)=|A|c$ for all $A\in\mathcal{L}_X^{e,e^*}$. If $A\in\mathcal{L}_X^e-\mathcal{L}_X^{e,e^*}$, we have $f_A(0)=0$. Indeed let $B\in\mathcal{L}_X^{e,e^*}$; using (b) we see that $|B|f_A(0)=0$ so that $f_A(0)=0$. We see that

$$f(x,e_0+x)=(-1)^{e^*(e)}c|A| \text{ if } x\in X_A, A\in\mathcal{L}_X^{e,e^*}), e_0\in E,$$

$$f(x,e_0+x)=0 \text{ if } x\in X_A, A\in\mathcal{L}_X^e-\mathcal{L}_X^{e,e^*}, e_0\in E.$$

This proves (a). The theorem is proved.

## 6. Simple $\mathcal{F}$-modules

**6.1.** For $(e,e^*)\in(E\times E^*)_X$ we denote by $((e,e^*))$ a simple module over the algebra $\bar{\mathcal{F}}^{e,e^*}$. Note that $((e,e^*))$ is well defined up to isomorphism and that

(a) $$\dim((e,e^*))=\sum_{A\in\mathcal{L}_X^{e,e^*}}|\bar{X}_A|,$$

see 4.6. Now $((e,e^*))$ can be regarded as a simple $\bar{\mathcal{F}}^e$-module on which $\bar{\mathcal{F}}^{e,e'^*}$ (with $e'^*\in E-\{e^*\}$) acts as zero and also as a simple $\mathcal{F}$-module via the isomorphism 4.2(a). (The summands $\bar{\mathcal{F}}^{e'},e'\ne e$ act on it as zero.) If $\mathcal{O}\in\overline{X^2}_{AB}$ and $\epsilon\in(A\cap B)^*$ then setting $\tilde{t}=\mathrm{tr}([\mathcal{O}]^\epsilon,((e,e^*)))$, we have:

(b) $$\tilde{t}=(-1)^{\epsilon(e)+e^*(e_0)}$$

if $A=B\in\mathcal{L}_X^{e,e^*}$ and $\mathcal{O}=\{(e_0+x,x),x\in o\}$ for some $e_0\in E$ and some $E$-orbit $o$ in $X_A=X_B$ and $\tilde{t}=0$ in all other cases.
(We use 4.5 and 4.1.) We see that

$$\{((e,e^*));(e,e^*)\in E\times E^*,e^*(e)=0,\mathcal{L}_X^{e,e^*}\ne\emptyset\}$$

is a complete set of mutually nonisomorphic simple $\mathcal{F}$-modules.

**Theorem 6.2.** *Let $M\subset P$ be subspaces of $E$ such that for some $D\in\mathcal{L}_X$ we have $M\subset D\subset E$. Let $P^\perp=\{e^*\in E^*;e^*|_P=0\}$. Then the $\mathcal{F}$-module $(MP)$ (see 4.7) is isomorphic to*

$$\oplus_{(e,e^*)\in M\times P^\perp}((e,e^*)).$$

*(Note that for $(e,e^*)\in M\times P^\perp$ we have automatically $(e,e^*)\in(E\times E^*)_X$ by the definition of $M,P$, so that $((e,e^*))$ is defined.) In particular, this $\mathcal{F}$-module is a direct sum of $|M||P^\perp|$ mutually nonisomorphic simple $\mathcal{F}$-modules.*

It is enough to show that for any $\mathcal{O}\in\overline{X^2}_{AB}$ and any $\epsilon\in(A\cap B)^*$, the sum

$$[t]=\sum_{(e,e^*)\in M\times E^*;e^*|_P=0}\mathrm{tr}([\mathcal{O}]^\epsilon,((e,e^*)))$$

is equal to $\mathrm{tr}([\mathcal{O}]^\epsilon,(MP))$.

Assume first that

(a) $A=B$ and $\mathcal{O}=\{(e_0+x,x),x\in o\}$ for some $e_0\in E$ and some $E$-orbit $o$ in $X_A=X_B$.

Using 6.1 we see that

$$[t]=\sum_{(e,e^*)\in M\times E^*;e\in B,e^*|_P=0,e^*|_B=0}(-1)^{\epsilon(e)+e^*(e_0)}$$
$$=\sum_{e\in M\cap B}(-1)^{\epsilon(e)})(\sum_{e^*\in E;e^*|_{B+P}=0}(-1)^{e^*(e_0)})$$
$$=\delta_{\epsilon|_{M\cap B}=0}|M\cap B||\delta_{e_0\in B+T}|E/(B+P)|.$$

Thus we have

$$[t]=|M\cap B||E|/|B+P|$$

if $\epsilon|_{M\cap B}=0$ and $e_0\in B+T$;

$$[t]=0$$

if $\epsilon|_{M\cap B}\neq 0$ or $e_0\notin B+T$. Next we assume that (a) is not satisfied. Then using 5.4 we see that $[t]=0$. It remains to use the formula for $\mathrm{tr}([\mathcal{O}]^\epsilon,(MP))$ given by 3.2.

**Colollary 6.3.** *In the setup of 6.2, the $\mathcal{F}$-module $\oplus_{(e,e^*)\in M\times P^\perp}((e,e^*))$ is positive.*

**6.4.** In the setup of 6.2 we have

(a) $$\dim(MP)=\sum_{A\in\mathcal{L}}|\bar{X}_A||A\cap M||E|/|A+P|.$$

Indeed, using 6.2 we have

$$\dim(MP)=\sum_{(e,e^*)\in M\times E^*;e^*|_P=0}\dim((e,e^*))).$$

Using now 6.1(a) we obtain

$$\begin{aligned}\dim(MP)&=\sum_{(e,e^*)\in M\times E^*;e^*|_P=0}\ \sum_{A\in\mathcal{L};e\in A,e^*|_A=0}|\bar{X}_A|\\&=\sum_{A\in\mathcal{L}}\ \sum_{(e,e^*);e\in M\cap A,e^*|_{P+A}=0}|\bar{X}_A|\\&=\sum_{A\in\mathcal{L}}|\bar{X}_A||M\cap A||E|/|P+A|,\end{aligned}$$

as required.

**6.5.** Taking $M=0,P=E$ in 6.2, 6.4 (with $X\neq\emptyset$), we have that $(0E)$ is a simple $\mathcal{F}$-module of dimension $\sum_{A\in\mathcal{L}}|\bar{X}_A|$.

**6.6.** Let $M'\subset M\subset P\subset P'$ be subspaces of $E$ such that for some $D\in\mathcal{L}_X$ we have $M\subset D\subset E$. Then the $\mathcal{F}$-modules $(MP)$, $(M'P')$ are defined. We show:

(a) The $\mathcal{F}$-module $(MP)$ contains a unique $\mathcal{F}$-submodule isomorphic to $(M'P')$. The existence of such a aubmodule follows from 6.2 since $e\in M'$ implis $e\in M$ and $e^*\in P'^\perp$ implies $e^*\in P^\perp$. The uniqueness follows from the fact that $(M,P)$ is multiplicity free.

**6.7.** Let $(e,e^*)\in(E\times E^*)_X$. Let $M_e$ be $\{0,e\}$ if $e\ne0$ and $M_e=0$ if $e=0$. We set $P_{e^*}=\ker\epsilon^*$. We have $M_e\subset P_{e^*}$ since $\epsilon^*(e)=0$. Let $A\in\mathcal{L}_X$ be such that $e\in A, e^*|A=0$. We have $M_e\subset A\subset P_{e^*}$ hence $(M_eP_{e^*})$ is defined. From 6.2 we have (with $\sim$ meaning "isomorphic to"):

$(M_eP_{e^*})\sim((e,e^*))\oplus((0,e^*))\oplus((e,0))\oplus((0,0))$ if $e\ne0, e^*\ne0$;
$(M_eP_0)\sim((e,0))\oplus((0,0))$ if $e\ne0$;
$(M_0P_{e^*})\sim((0,e^*))\oplus((0,0))$ if $e^*\ne0$;
$(M_0P_0)\sim((0,0))$.

It follows that

$((e,e^*))\sim(M_eP_{e^*})/)(M_0P_{e^*})+(M_eP_0))$ if $e\ne0, e^*\ne0$;
$((e,0))\sim(M_eP_0)/(M_0P_0)$ if $e\ne0$;
$((0,e^*))\sim(M_0P_{e^*})/(M_0P_0)$ if $e^*\ne0$;
$((0,0))\sim(M_0P_0)$.

## 7. Positivity for Weyl group modules

**7.1.** Let $W$ be an irreducible Weyl group and let $J$ be the corresponding asymptotic Iwahori-Hecke ring (see [L87, §2]); let $\mathbf{J}=\mathbf{C}\otimes J$, an algebra over $\mathbf{C}$ with a basis $\{t_w;w\in W\}$. We have a direct sum decomposition $\mathbf{J}=\oplus_{c\in ce(W)}\mathbf{J}_c$ (as algebras) where $ce(W)$ is the set of two-sided cells of $W$ and $\mathbf{J}_c$ is spanned by $\{t_w;w\in c\}$. It is known that there is a canonical algebra isomorphism $\mathbf{C}[W]\xrightarrow{\sim}\mathbf{J}$, see [L87, 3.1]. Via this isomorphism any $\mathbf{C}[W]$-module $\mathcal{E}$ can be regarded as a $\mathbf{J}$-module $\mathcal{E}_\infty$. For $c\in ce(W)$ let $Mod_c(W)$ be the category of $\mathbf{C}[W]$-modules such that $\mathbf{J}_{c'}\mathcal{E}_\infty=0$ for any $c'\in ce(W)-\{c\}$. Note that any irreducible $\mathbf{C}[W]$-module is in $Mod_c(W)$ for a unique $c\in ce(W)$.

**7.2.** Let $c\in ce(W)$ and let $\mathcal{E}\in Mod_c(W)$. We say that $\mathcal{E}$ is *positive* if there exists a $\mathbf{C}$-basis of $\mathcal{E}_\infty$ such that any $t_w (w\in c)$ acts in this basis through a matrix all of whose entries are in $\mathbf{R}_{\ge0}$. In [L18] it is shown that if $\mathcal{E}\in Mod_c(W)$ is assumed to be irreducible then $\mathcal{E}$ is positive if and only if $\mathcal{E}$ is special. (In the case where $W$ is of classical type, the fact that the special representations of $W$ are positive was proved earlier, see [L16].)

**7.3.** Now let $c\in ce(W)$. In [L87a] it was conjectured (and partially proved) that

(a) there exists an algebra isomorphism $\mathbf{J}_c\to\mathbf{C}\otimes K_{\mathcal{G}}(X^2)$ (where $\mathcal{G}$ is a finite group attached to $c$ and $X$ is a certain explicit finite set with $\mathcal{G}$ action) carrying the basis 0.1(a) of $\underline{K}_{\mathcal{G}}(X^2)$ onto the basis $\{t_w;w\in c\}$ of $\mathbf{J}_c$.
(The proof of (a) was completed in [BFO].)

For $\mathcal{E}\in Mod_c(W)$ we denote by $\mathcal{E}^!_\infty$ the $\underline{K}_{\mathcal{G}}(X^2)$-module corresponding to the $\mathbf{J}_c$-module $\mathcal{E}_\infty$ under the isomorphism (a).

**7.4.** In [L19] a set of (not necessarily irreducible) representations in $Mod_c(W)$ (called the new representations) was introduced; these new representations form a $\mathbf{Z}$-basis of the Grothendieck group of $Mod_c(W)$. For example the special representations and the constructible representations of $W$ in $Mod_c(W)$ are new, but in general they are not the only new ones.

**7.5.** We now assume that $W$ is of type $B, C$ or $D$ so that $\mathcal{G}$ in 7.3 is $E$ in 1.1. It is known that the isomorphism 7.3(a) has the following property.

(a) We have $\mathcal{L}_X = \mathcal{L}'(E)$ where $\mathcal{L}'(E)$ is the set of subgroups of $E$ attached in [L87a] to the various left cells contained in $c$.

Let $\mathfrak{F}_c(W)$ be the set of all $\mathcal{E} \in Mod_c(W)$ such that $\mathcal{E}^!_\infty = (MP)$ with $M, P$ as in 4.7. From 4.7(a) it follows that

(b) if $\mathcal{E} \in \mathfrak{F}_c(W)$ then $\mathcal{E}$ is positive.

If $\mathcal{E} \in \mathrm{Mod}_c(W)$ is irreducible then $\mathcal{E}^!_\infty = ((e, e^*))$ for a well defined $(e, e^*) \in (E \times E^*)_X$. In [L19] to $\mathcal{E}$ we have attached two subgroups $M, P$ in $\mathcal{L}'(E) = \mathcal{L}_X$ such that $M \subset P$ and $(e, e^*) \in M \oplus P^\perp \subset (E \times E^*)_X$. Define $\tilde{\mathcal{E}} \in Mod_c(W)$ by $\tilde{\mathcal{E}}^!_\infty = \oplus_{(e,e^*) \in M \oplus P^\perp} ((e, e^*))$. Then $\tilde{\mathcal{E}}$ is a new representation in $Mod_c(W)$ in the sense of [L19] and all new representations in $Mod_c(W)$ are obtained in this way. From 6.3 or (b) we see that $\tilde{\mathcal{E}}$ is positive. Thus

(c) all new representations in $Mod_c(W)$ are positive.

For $W$ of type $A$, (c) is also true. We expect that (c) remains true for $W$ of exceptional type.

Department of Mathematics, M.I.T., Cambridge, MA 02139